\documentclass{article}

\usepackage{tikz,graphics,xcolor}
\usepackage{amsmath,amsthm,mathtools,amsfonts,amssymb}
\usepackage{enumitem,comment}
\allowdisplaybreaks

\newtheorem{theorem}{Theorem}
\newtheorem*{theorem*}{Theorem}
\newtheorem{proposition-definition}{Proposition-Definition}
\newtheorem{definition}[theorem]{Definition}

\newtheorem{lemma}[theorem]{Lemma}

\theoremstyle{remark}
\newtheorem{remark}[theorem]{Remark}
\newcommand{\R}{\mathbb{R}}
\newcommand{\Qp}{\mathbb{Q}_p}

\newcommand{\Z}{\mathbb{Z}}
\newcommand{\Q}{\mathbb{Q}}
\newcommand{\acm}{\overline{\mathrm{ac}}_m}

\newcommand{\ord}{\mathrm{ord}}

\renewcommand{\d}{\mathrm{d}}
\newcommand{\id}{\mathrm{id}}

\newcommand{\abs}[1]{\lvert #1\rvert}
\newcommand{\restr}[2]{\ensuremath{#1\big|_{#2}}}
\newcommand{\deriv}[2]{\ensuremath{d#1/d#2}}

\title{Lipschitz extensions of definable $p$-adic functions}
\author{Tristan Kuijpers}
\date{}
\begin{document}
\maketitle
\abstract{In this paper, we prove a definable version of Kirszbraun's theorem in a non-Archimedean setting for definable families of functions in one variable. More precisely, we prove that every definable function $f:X\times Y\to\Q_p^s$, where $X\subset \Q_p$ and $Y\subset \Q_p^r$, that is $\lambda$-Lipschitz in the first variable, extends to a definable function $\tilde{f}:\Q_p\times Y\to\Q_p^s$ that is $\lambda$-Lipschitz in the first variable.}
\section{Introduction}
In 1934, Kirszbraun proved that every $\lambda$-Lipschitz function $f:S\subset \R^r\to\R^s$ extends to a $\lambda$-Lipschitz function $\tilde{f}:\R^r\to\R^s$ (see \cite{kirszbraun}). In 1983, Bhaskaran proved that a version of Kirszbraun's theorem still holds in a non-Archimedean setting, more precisely, for all spherically complete fields (see \cite{bhaskaran}). Recently, in 2010, Aschenbrenner and Fischer proved a definable version of Kirszbraun's theorem. In particular, they proved that every $\lambda$-Lipschitz function $f:S\subset\R^r\to\R^s$, that is definable in an expansion of the ordered field of real numbers, extends to a $\lambda$-Lipschitz function $\tilde{f}:\R^r\to\R^s$ that is definable in the same structure (see \cite{aschenbrenner}).

The proof of Bhaskaran relies in an essential way on Zorn's Lemma, which makes it far from being applicable to a definable setting. Therefore Aschenbrenner posed the question whether there could be a \emph{definable} version of Kirszbraun's theorem in a non-Archimedean setting. In this paper we partially answer that question and prove a definable version of Kirszbraun's theorem in a non-Archimedean setting for definable families of functions in one variable. More precisely, we prove that every definable function $f:X\times Y\to\Q_p^s$, where $X\subset \Q_p$ and $Y\subset \Q_p^r$, that is $\lambda$-Lipschitz in the first variable, extends to a definable function $\tilde{f}:\Q_p\times Y\to\Q_p^s$ that is $\lambda$-Lipschitz in the first variable. By \emph{definable}, we mean definable in either a semi-algebraic or a subanalytic structure on $\Q_p$. Working with these languages will allow us to use a cell decomposition result (see Theorem \ref{thm_preparation}) that is essential for the construction of Lipschitz extensions.

In a first approach we use a more easy construction to obtain a $\Lambda$-Lipschitz extension, where $\Lambda$ is possibly larger than $\lambda$. In a second and more involved approach, we show one can take $\Lambda$ equal to $\lambda$. More generally, we prove our results for finite field extensions of $\Qp$.

\subsection*{Acknowledgments}
The autor would like to thank Raf Cluckers for proposing the idea for this paper, the many fruitful discussions and the constant optimism during the preparation of this paper. The author is also grateful to Matthias Aschenbrenner, for it was his question that formed the inspiration for this research project.

\section{Preliminary definitions and facts}
Let $p$ be a prime number, and $\Q_p$ the field of $p$-adic numbers. Let $K$ be a finite field extension of $\Q_p$. Denote by $\ord:K^\times \to \Z$ the valuation. Denote by $\mathcal{O}_K$ the valuation ring, by $\mathcal{M}_K$ the maximal ideal of $\mathcal{O}_K$ and by $\pi_K$ a fixed generator of $\mathcal{M}_K$. Let $q$ denote the cardinality of the residue field. Finally, let $\acm: K\to \mathcal{O}_K/(\pi_K^m)$ be the angular component map of depth $m$, sending every nonzero $x$ to $x\pi_K^{-\ord(x)} \mod (\pi_K^m)$ and 0 to 0.

The valuation induces a non-Archimedean norm on $K$ by setting $\abs{x}=q^{-\ord(x)}$ for nonzero $x$, and $\abs{0}=0$. This extends to a non-Archimedean norm on $K^s$ by setting $\abs{(x_1,\ldots,x_s)} = \max_i \{\abs{x_i}\}$. A function $f:K^r\to K^s$ is said to be $\lambda$-Lipschitz, with $\lambda\in\R$, if $\abs{f(x)-f(y)}\leq \lambda \abs{x-y}$ for all $x,y\in K^r$. One calls $\lambda$ the \emph{Lipschitz constant} of $f$.

Say a set $X\subset K^r$ is \emph{definable} if it is definable in either a semi-algebraic or a subanalytic structure on $K$. This means that $X$ is given by a first-order formula, possibly with parameters form $K$, in the semi-algebraic or subanalytic language (see \cite{ccl} for more details). For the convenience of the reader, we recall these languages. The semi-algebraic (or \emph{Macintyre}) language is the language $\mathcal{L}_{\text{Mac}}=(+,-,\cdot,\{P_n\}_{n>0},0,1)$, where the predicates $P_n$ stand for the $n$-th powers in $K$. The subanalytic language is the language $\mathcal{L}_{\text{an}}=\mathcal{L}_{\text{Mac}} \cup (^{-1},\cup_{m>0}K\{x_1,\ldots,x_m\})$, where $^{-1}$ is interpreted as the multiplicative inverse extended by $0^{-1}=0$, and where every function symbol from $K\{x_1,\ldots,x_m\}$ is interpreted as the restricted analytic function $K^m\to K$ given by
\[x\mapsto \begin{cases} f(x) &\text{if }x\in \mathcal{O}_K^m,\\0&\text{otherwise},\end{cases}\]
where $f$ is a formal power series converging on $\mathcal{O}_K^m$. Let $X\subset K^r$ be a definable set, then a function $f:X\to K^s$ is definable if its graph is a definable subset of $K^{r+s}$.

We work with the notion of $p$-adic cells as given in \cite{ch}. We recall the main definitions and properties. For integers $m,n>0$, let $Q_{m,n}$ be the (definable) set
\[Q_{m,n} = \{x\in K^\times \mid \ord(x)\in n\Z,\ \acm(x)=1\}.\]
\begin{definition}
Let $Y$ be a definable set. A \emph{cell} $C\subset K\times Y$ over $Y$ is a (nonempty) set of the form
\begin{equation*}
	\resizebox{.85\hsize}{!}{$C = \{(x,y)\in K\times Y\mid y\in Y',\ \abs{\alpha(y)}\mathrel{\square_1}\abs{x-c(y)}\mathrel{\square_2} \abs{\beta(y)},\ x-c(y)\in\xi Q_{m,n}\}$,}
\end{equation*}
where $Y'\subset Y$ is a definable set, $\xi\in K$, $\alpha, \beta: Y'\to K^\times$ and $c:Y'\to K$ are definable functions, $\square_i$ is either $<$ or ``no condition'', and such that $C$ projects surjectively onto $Y'$. We call $c$ and $\xi Q_{m,n}$ the \emph{center} and the \emph{coset} of the cell $C$, respectively. If $\xi = 0$ we call $C$ a \emph{0-cell}, otherwise we call $C$ a \emph{1-cell}. We call $Y'$ the \emph{base} of the cell $C$.
\end{definition}
\begin{definition}
	Let $Y$ be a definable set. Let $C\subset K\times Y$ be a 1-cell over $Y$ with center $c$ and coset $\xi Q_{m,n}$. Then, for each $(t,y) \in C$ with $y\in Y$, there exists a unique maximal ball $B$ containing $t$ and satisfying $B\times \{y\}\subset C$, where maximality is under inclusion. If $\ord(t-c(y))=l$, this ball is of the form
		\[B = B_{l,c(y),m,\xi} = \{x\in K\mid \ord(x-c(y)) = l,\, \acm(x-c(y)) = \acm(\xi)\}.\]

	We call the collection of all these maximal balls the \emph{balls of the cell $C$}. For fixed $y_0\in Y$, we call the collection of balls $\{B_{l,c(y_0),m,\xi}\mid B_{l,c(y_0),m,\xi}\times \{y_0\}\subset C\}$ the balls of the cell $C$ \emph{above $y_0$}. If $C\subset K\times Y$ is a $0$-cell, we define the collection of balls of $C$ to be the empty collection.
\end{definition}
Notice that $B_{l,c(y),m,\xi}$ is a ball of diameter $q^{-(l+m)}$, in particular, for every $x_1,x_2\in B_{l,c(y),m,\xi}$ it holds that $\abs{x_1-x_2}\leq q^{-(l+m)}$.
\begin{definition}[Jacobian property]
	Let $f:B\to B'$ be a function, where $B,B'\subset K$ are balls. Say that $f$ has the \emph{Jacobian property} if the following conditions hold:
	\begin{enumerate}
		\item $f$ is a bijection;
		\item $f$ is continuously differentiable on $B$, with derivative $\deriv{f}{x}$;
		\item $\ord (\deriv{f}{x})$ is constant (and finite) on $B$;
		\item for all $x,y\in B$ with $x\neq y$, one has:
			\[\ord(f(x)-f(y)) = \ord(\deriv{f}{x}) + \ord(x-y).\]
	\end{enumerate}
\end{definition}
\begin{definition}
	Let $f:S\subset K\times Y\to K$ be a function. Then we define
	\[f\times \id: S\to K\times Y: (x,y)\mapsto (f(x,y),y),\]
	and we denote with $S_f$ the image of $f\times\id$.
\end{definition}
\begin{definition}
	Let $f:S\subset K\times Y\to K^s$ be a function. Then we define for every $y\in Y$
		\[f_y: S_y\to K^s: x\mapsto f(x,y),\]
	where $S_y$ denotes the fiber $S_y = \{x\in K\mid (x,y)\in S\}$. 
\end{definition}
\begin{definition}
Let $Y$ be a definable set, let $C\subset K\times Y$ be a 1-cell over $Y$, and let $f:C\to K$ be a definable function. Say that $f$ is \emph{compatible} with the cell $C$ if either $C_f$ is a 0-cell over $Y$, or the following holds: $C_f$ is a 1-cell over $Y$ and for each $y\in Y$ and each ball $B$ of $C$ above $y$ and each ball $B'$ of $C_f$ above $y$, the functions $\restr{f_y}{B}$ and $\restr{f_y^{-1}}{B'}$ have the Jacobian property.

If $g:C\to K$ is a second definable function which is compatible with the cell $C$ and if we have $C_f = C_g$ and $\ord(\frac{\partial f(x,y)}{\partial{x}})= \ord (\frac{\partial g(x,y)}{\partial x})$ for every $(x,y)\in C$, then we say that $f$ and $g$ are \emph{equicompatible} with $C$.

If $C'\subset K\times Y$ is a 0-cell over $Y$, any definable function $h:C'\to K$ is said to be compatible with $C'$, and $h$ and $k:C'\to K$ are equicompatible with $C'$ if and only if $h=k$.
\end{definition}
The following theorem is based on Theorem 3.3 of \cite{ch}. This theorem is the result of a constant refinement of the concept of $p$-adic cell decomposition for semi-algebraic and subanalytic structures. Earlier versions are due to Cohen \cite{cohen}, Denef \cite{denef-84,denef}, Cluckers \cite{cluckers}, and relate to the quantifier elimination results from Macintyre \cite{mac} and  Denef-van den Dries \cite{denef-vdd}. 
\begin{theorem}\label{thm_preparation}
Let $S\subset K\times Y$ and $f:S\to K$ be definable. Then there exists a finite partition of $S$ into cells $C$ over $Y$ such that the restriction $\restr{f}{C}$ is compatible with $C$ for each cell $C$. Moreover, for each cell $C$ there exists a definable function $m:C\to K$, a definable function $e:Y\to K$ and coprime integers $a$ and $b$ with $b>0$, such that for all $(x,y)\in C$
\[m(x,y)^b = e(y)(x-c(y))^a,\]
where $c$ is the center of $C$, and such that if one writes $c'$ for the center of $C_f$, one has that $g = m+c'$ and $f$ are equicompatible with $C$ (we use the conventions that $b=1$ whenever $a=0$, that $a=0$ whenever $C$ is a 0-cell, and that $0^0=1$). 

Furthermore, if $C$ and $C_f$ are 1-cells, then for every $y\in Y$ one has that $f_y(B) = g_y(B)$ for every ball $B$ of $C$ above $y$, and the formula
\begin{equation}\label{eq_formula_order}
\ord\left(\frac{\partial f(x,y)}{\partial{x}}\right) = \ord(e(y)^{1/b}q) + (q-1)\ord(x-c(y))
\end{equation}
holds for all $(x,y)\in C$, where $q=a/b$ and where we use the convenient notation $\ord(t^{1/b}) = \ord(t)/b$, for $t\in K$ and $b>0$ a positive integer.
\end{theorem}

%
\begin{proof}
	The existence of a finite partition of $S$ in cells $C$ over $Y$, and for every such a cell $C$ the existence of $g=m+c'$ such that $f$ and $g$ are equicompatible with $C$, follows immediately from Theorem 3.3 in \cite{ch}. 
	
	Now assume that $C$ and $C_f$ are 1-cells. It is easy to see that $f_y(B) = g_y(B)$ for every $y\in Y$ and every ball $B$ of $C$ above $y$. 
	
	We now prove \eqref{eq_formula_order}. Fix $(x,y)\in C$. Since $f$ and $g$ are equicompatible, we have $\ord(\frac{\partial f(x,y)}{\partial{x}}) = \ord(\frac{\partial g(x,y)}{\partial{x}})$, so we only need to prove that \eqref{eq_formula_order} holds for $f$ replaced by $g$. For this, we first note that
	\begin{equation}\label{eq_ord_g}
	\ord(g(x,y)-c'(y)) = [\ord(e(y))+a\cdot\ord(x-c(y))]/b.
	\end{equation}
	It is also immediate that
	\begin{equation}\label{eq_g_rechts}
	\ord\left(\frac{\partial ([g(x,y)-c'(y)]^b)}{\partial{x}}\right) = \ord(e(y)a(x-c(y))^{a-1}).
	\end{equation}
	On the other hand, by the chain rule, the left hand side of \eqref{eq_g_rechts} also equals
	\begin{equation}\label{eq_g_links}
	\ord\left(\frac{\partial ([g(x,y)-c'(y)]^b)}{\partial{x}}\right) =\ord\left(b[g(x,y)-c'(y)]^{b-1}\frac{\partial g(x,y)}{\partial{x}}\right).
	\end{equation}
	Equating the right hand sides of \eqref{eq_g_rechts} and \eqref{eq_g_links}, and using \eqref{eq_ord_g}, one easily finds the required formula.
\end{proof}

Comparing sizes of balls between which there is a function with the Jacobian property, we obtain the following useful formula. 
\begin{lemma}\label{lemma_l'm'ordf'lm}
	Let $f:B_{l,c(y),m,\xi}\to B_{l',c'(y),m',\xi'}$ be a function with the Jacobian property. Then $l'+m' = \ord(df/dx)+l+m$.
\end{lemma}

\section{Existence of Lipschitz extensions}
We now proceed towards proving the existence of definable Lipschitz extensions of definable families of functions in one variable. Let us first formulate the main theorem of this paper:
\begin{theorem}\label{thm_main_param_const1}
	Let $Y\subset K^r$ and $X\subset K$ be definable sets and let $f:X\times Y\to K^s$ be a definable function that is $\lambda$-Lipschitz in the first variable. Then $f$ extends to a definable function $\tilde{f}:K\times Y\to K^s$ that is $\lambda$-Lipschitz in the first variable, i.e. $\tilde{f}_y$ is $\lambda$-Lipschitz for every $y\in Y$.
\end{theorem}
\begin{remark}\label{remark_lipschitz_constant}
By rescaling, if suffices to proof the theorem for $\lambda=1$. Also, since we use the max-norm on $K^s$, it is enough to prove the theorem for $s=1$.
\end{remark}

Firstly, we present a very general way of \emph{gluing} Lipschitz extensions of a given function to obtain a Lipschitz extension with a larger domain (this is Lemma \ref{lemma_glue}).

Secondly, given a definable function that is $\lambda$-Lipschitz in the first variable, we give a more easy construction to obtain a definable extension that is $\Lambda$-Lipschitz in the first variable, where $\Lambda$ is possibly larger than $\lambda$ (this is Theorem \ref{thm_main_param}).

Thirdly and lastly, using a more involved argument, we show that one can take $\Lambda$ equal to $\lambda$ (this is Theorem \ref{thm_main_param_const1}).
\begin{lemma}[Gluing extensions]\label{lemma_glue}
Let $X\subset K^r$ be a definable set and let $f:X\to K$ be a definable and $\lambda$-Lipschitz function. Let $X = \cup_{i=1}^k X_i$ be a finite covering of $X$ by definable subsets $X_i$. Call $f_i = \restr{f}{X_i}: X_i\to K$.

If every $f_i$ extends to a definable and $\Lambda_i$-Lipschitz map $\tilde{f_i}:K^r\to K$, with $\Lambda_i\geq\lambda$, then $f$ extends to a definable and $\Lambda$-Lipschitz map $\tilde{f}:K^r\to K$, where $\Lambda = \max_i \{\Lambda_i\}$.
\end{lemma}
\begin{proof}
We prove the lemma first for $k=2$. Define $T_1 = \{ x\in K^r\mid \d(x,X_1)\leq \d(x,X_2)\}$, where $\d(x,A)$ denotes the distance from $x$ to the set $A$, i.e. $\d(x,A) = \inf\{\abs{x-a}\mid a\in A\}$. Define $T_2 = K^r\setminus T_1$, and let
\[ \tilde{f}:K^r\to K: x\mapsto \begin{cases} \tilde{f_1}(x)&\text{if }x\in T_1,\\\tilde{f_2}(x)&\text{if }x\in T_2.\end{cases}\]
	Clearly $\tilde{f}$ is a definable extension of $f$. We prove that $\tilde{f}$ is $\Lambda$-Lipschitz, where $\Lambda = \max\{\Lambda_1,\Lambda_2\}$. The only nontrivial fact to verify is that for $t_1\in T_1$ and $t_2\in T_2$, we have $\abs{\tilde{f}(t_1)-\tilde{f}(t_2)}\leq \Lambda\abs{t_1-t_2}$.

Since every definable and $\lambda$-Lipschitz function extends uniquely to a definable and $\lambda$-Lipschitz function on the topological closure of its domain, we may assume that $X$, $X_1$ and $X_2$ are topologically closed. 


Fix elements $a_i\in X_i$ such that $\abs{t_i-a_i} = \d(t_i,X_i)$, for $i=1,2$. 
It then always holds that
	\begin{equation}\label{eq_a2x2<a1a2}
		\abs{t_2-a_2}<\abs{a_1-a_2}.
	\end{equation}
We can now calculate as follows: 
\begin{align*}
	\abs{\tilde{f}(t_1)-\tilde{f}(t_2)} &= \abs{\tilde{f_1}(t_1)-f(a_1)+f(a_2)-\tilde{f_2}(t_2)+f(a_1)-f(a_2)}\\
	&\leq\max\{\abs{\tilde{f_1}(t_1)-f(a_1)},\abs{f(a_2)-\tilde{f_2}(t_2)},\abs{f(a_1)-f(a_2)}\}\\
	&\leq\max\{\Lambda_1\abs{t_1-a_1},\Lambda_2\abs{a_2-t_2},\lambda\abs{a_1-a_2}\}\\
	&\leq\max\{\Lambda_1,\Lambda_2\}\max\{\abs{t_1-a_1},\abs{a_2-t_2},\abs{a_1-a_2}\}\\
	&\stackrel{\eqref{eq_a2x2<a1a2}}{=}\Lambda\max\{\abs{t_1-a_1},\abs{a_1-a_2}\}\\
	&\leq\Lambda\abs{t_1-t_2},
\end{align*}
where we only have to verify the last inequality. For this we prove that 
\begin{equation}\label{eq_preparatory}\max\{\abs{t_1-a_1},\abs{a_1-a_2}\}\leq \abs{t_1-t_2},\end{equation}
by considering two cases.
\begin{description}
\item[Case 1: $\abs{t_1-a_1}<\abs{a_1-a_2}$.] 
	It then holds that
	\begin{equation}\label{eq_a1a2=x1a2N}
		\abs{a_1-a_2}=\abs{t_1-a_2}.
	\end{equation}
	So
	\begin{equation}\label{eq_lange}
		\abs{t_2-a_2}\stackrel{\eqref{eq_a2x2<a1a2}}{<}\abs{a_1-a_2}\stackrel{\eqref{eq_a1a2=x1a2N}}{=}\abs{t_1-a_2},
	\end{equation}
	hence
	\begin{equation}
		\abs{a_1-a_2}\stackrel{\eqref{eq_a1a2=x1a2N}}{=}\abs{t_1-a_2}\stackrel{\eqref{eq_lange}}{=}\abs{t_1-t_2}.
	\end{equation}	
\item[Case 2: $\abs{a_1-a_2}\leq\abs{t_1-a_1}$.] Suppose that
\begin{equation}\label{eq_x1x2<x1a1} \abs{t_1-t_2}<\abs{t_1-a_1},\end{equation}
then 
\begin{equation}\label{eq_x1a1=x2a1=a1a2}\abs{t_1-a_1}\stackrel{\eqref{eq_x1x2<x1a1}}{=}\abs{t_2-a_1} \stackrel{\eqref{eq_a2x2<a1a2}}{=} \abs{a_1-a_2}.\end{equation}
By the choice of $a_1$ and the fact that $t_1\in T_1$, we know $\abs{t_1-a_2}\geq \abs{t_1-a_1}$, so by \eqref{eq_x1a1=x2a1=a1a2} equality holds:
\begin{equation}\label{eq_x1a1=x1a2}\abs{t_1-a_1}=\abs{t_1-a_2}.\end{equation}
Together with \eqref{eq_x1x2<x1a1}, this implies
\begin{equation}\label{eq_x1x2<x1a2twee}\abs{t_1-t_2}<\abs{t_1-a_2}.\end{equation}
So finally,
\begin{equation*}
	\abs{a_1-a_2}\stackrel{\eqref{eq_x1a1=x2a1=a1a2}}{=}\abs{t_1-a_1}\stackrel{\eqref{eq_x1a1=x1a2}}{=}\abs{t_1-a_2}\stackrel{\eqref{eq_x1x2<x1a2twee}}{=}\abs{t_2-a_2} \stackrel{\eqref{eq_a2x2<a1a2}}{<}\abs{a_1-a_2},
\end{equation*}
which is a contradiction.
\end{description}
This proves \eqref{eq_preparatory}, and therefore the lemma is proved for $k=2$. An easy induction argument then proves the lemma for general $k$.
\end{proof}
\begin{remark} Lemma \ref{lemma_glue} remains true is one replaces every instance of the word ``Lipschitz'' by ``Lipschitz in the first variable''.
\end{remark}

\begin{theorem}\label{thm_main_param}
	Let $Y\subset K^r$ and $X\subset K$ be definable sets and let $f:S=X\times Y\to K^s$ be a definable function that is $\lambda$-Lipschitz in the first variable. Then there exists $\Lambda\geq \lambda$ such that $f$ extends to a definable function $\tilde{f}:K\times Y\to K^s$ that is $\Lambda$-Lipschitz in the first variable, i.e. $\tilde{f}_y$ is $\Lambda$-Lipschitz for every $y\in Y$.
\end{theorem}
\begin{proof}
	By Remark \ref{remark_lipschitz_constant}, we may assume that $\lambda=1$ and $s=1$. By Theorem \ref{thm_preparation} and (the remark after) Lemma \ref{lemma_glue}, we may assume that $S$ is a cell over $Y$ with which $f$ is compatible. Furthermore, we may assume that the base of $S$ is $Y$.

If $S_f$ is a $0$-cell over $Y$ with center $c'$, we define 
\[\tilde{f}:K\times Y\to K: (x,y)\mapsto c'(y).\]
Clearly, $\tilde{f}$ is a definable extension of $f$ and for all $y\in Y$, $\tilde{f}_y$ is $1$-Lipschitz.

Assume from now on that $S$ and $S_f$ are 1-cells over $Y$, with center $c$ and $c'$, and coset $\xi Q_{m,n}$ and $\xi' Q_{m',n'}$, respectively. 

	We define $\tilde{f}$ as follows:

	\[\tilde{f}:K\times Y\to K: (x,y)\mapsto\begin{cases} f(x,y) & \text{if } (x,y)\in S,\\c'(y)&\text{if }(x,y)\not\in S.\end{cases}\]
Clearly, $\tilde{f}$ is a definable extension of $f$. We prove that $\tilde{f}_y$ is $q^{m'}$-Lipschitz for every $y\in Y$.

Fix $y\in Y$. Let $t_1\in X$ and $t_2\not \in X$. Let $l$ and $l'$ be such that $t_1\in B_{l,c(y),m,\xi}$ and $f(t_1,y)\in B_{l',c'(y),m',\xi'}$. Then 
	\begin{align} \abs{f(t_1,y)-c'(y)} &= q^{-l'}\notag\\
			&= q^{-\ord(\partial f(t_1,y)/\partial x)}q^{m'-m}q^{-\ord(t_1-c(y))}\notag\\
			&\leq q^{m'-m}\abs{t_1-c(y)},\label{eq_mm'xc}
	\end{align}
	where the second equality follows from Lemma \ref{lemma_l'm'ordf'lm} and the last inequality holds because $f$ is $1$-Lipschitz in the first variable, and therefore $\abs{\partial f(t_1,y)/\partial x}\leq 1$.
	There are two cases to consider.
	\begin{description}
		\item[Case 1: $\abs{t_1-c(y)} = \abs{t_2-c(y)}$.] Because $B_{l,c(y),m,\xi}$ is a ball of diameter $q^{-m-l}$, it holds that $\abs{t_1-t_2}> q^{-m-l}$, or put differently:
			\begin{equation}\label{eq_l<mx1x2}
				q^{-m}\abs{t_1-c(y)}<\abs{t_1-t_2}.
			\end{equation}
			Therefore
			\begin{align*}
                        \abs{\tilde{f}_y(t_1)-\tilde{f}_y(t_2)} &= \abs{\tilde{f}(t_1,y)-\tilde{f}(t_2,y)}\\
                        &= \abs{f(t_1,y)-c'(y)}\\
                        &\stackrel{\eqref{eq_mm'xc}}{\leq} q^{m'-m}\abs{t_1-c(y)}\\
                        &\stackrel{\eqref{eq_l<mx1x2}}{<} q^{m'}\abs{t_1-t_2}.
			\end{align*}
		\item[Case 2: $\abs{t_1-c(y)} \neq \abs{t_2-c(y)}$.] From the non-Archimedean property it then follows that
			\begin{equation}\label{eq_x1c<x1x2}
				\abs{t_1-c(y)}\leq\abs{t_1-t_2},
			\end{equation}
			so we find
			\begin{align*}
                        \abs{\tilde{f}_y(t_1)-\tilde{f}_y(t_2)} &= \abs{\tilde{f}(t_1,y)-\tilde{f}(t_2,y)}\\
                        &=\abs{f(t_1,y)-c'(y)} \\
                        &\stackrel{\eqref{eq_mm'xc}}{\leq} q^{m'-m}\abs{t_1-c(y)}\\
                        &\stackrel{\eqref{eq_x1c<x1x2}}{\leq}q^{m'-m}\abs{t_1-t_1}.\qedhere
                        \end{align*}
	\end{description}
\end{proof}
\begin{remark}\label{remark_M}
	Analyzing the proof of Theorem \ref{thm_main_param}, we find that one can take $\Lambda = \lambda\max_i\{q^{m_i'}\}$, where $\lambda$ is the Lipschitz constant of $f$ (in the first variable), and the $m_i'$ correspond to the 1-cells in the cell decomposition of $S_f$.
\end{remark}
\begin{remark}\label{rem_phi}
	We can even improve (i.e. decrease) $\Lambda$ from Remark \ref{remark_M} as follows. In the proof, the worst Lipschitz constant occurs in \textbf{Case 1}. We can get around this case in the following way (as in the beginning of Theorem \ref{thm_main_param}, we assume that $S$ and $S_f$ are $1$-cells over $Y$ with center $c$ and coset $\xi Q_{m,n}$).

	For every nonzero $a\in \mathcal{O}_K^\times/(\pi_K^m)$, choose $\xi_m(a)\in \mathcal{O}_K^\times$ to be a class representative of $a$. Since we only need to make a finite number of representative choices, $\xi_m:\mathcal{O}_K^\times/(\pi_K^m)\to K$ is a definable map. Let $\varphi:K\times Y\to K$ be the definable map rescaling the angular component as follows:
	\begin{align*}
	\varphi:K\times Y&\to K:\\
	(x,y)&\mapsto \begin{cases}(x-c(y))\xi_m(\acm(x-c(y))^{-1}\acm(\xi))+c(y) & \text{if }x\neq c(y),\\ c(y) & \text{if }x=c(y).\end{cases}
	\end{align*}
	It is not difficult to see that for every $y\in Y$, $\varphi_y:K\to K$ is $1$-Lipschitz. Now let $\hat{f}$ be the extension described in the proof of Theorem \ref{thm_main_param} in the case that $S$ and $S_f$ are $1$-cells over $Y$ (remark that in Theorem \ref{thm_main_param}, this extension is denoted with $\tilde{f}$). Then $\tilde{f}:K\times Y\to K: (x,y)\mapsto \hat{f}(\varphi(x,y),y)$ is a definable extension of $f$ that is $q^{m'-m}$-Lipschitz in the first variable. One can therefore take $\Lambda = \lambda\max_i\{q^{m_i'-m_i}\}$, where $\lambda$ is the Lipschitz constant of $f$ (in the first variable), and the $m_i$ and $m_i'$ correspond to the 1-cells in the cell decomposition of $S$ and $S_f$, respectively.
\end{remark}

Note that in the proof of Theorem \ref{thm_main_param} we did not use the full generality of Theorem \ref{thm_preparation}. We will now prove Theorem \ref{thm_main_param_const1}, the main theorem of this paper, which uses a more involved extension for which the Lipschitz constant doesn't grow. For this, the full power of Theorem \ref{thm_preparation} is used. Again, the result is formulated in definable families of functions. For clarity, we repeat the formulation of Theorem \ref{thm_main_param_const1}. 

\begin{theorem*}
	Let $Y\subset K^r$ and $X\subset K$ be definable sets and let $S=X\times Y$. Let $f:S\to K^s$ be a definable function that is $\lambda$-Lipschitz in the first variable. Then $f$ extends to a definable function $\tilde{f}:K\times Y\to K^s$ that is $\lambda$-Lipschitz in the first variable, i.e. $\tilde{f}_y$ is $\lambda$-Lipschitz for every $y\in Y$.
\end{theorem*}
	\begin{proof}
		By Remark \ref{remark_lipschitz_constant}, we may assume that $\lambda=1$ and $s=1$. By Theorem \ref{thm_preparation} and (the remark after) Lemma \ref{lemma_glue}, we may assume that $S$ is a cell over $Y$ with which $f$ is compatible. Furthermore, we may assume that the base of $S$ is $Y$.
				
If $S_f$ is a $0$-cell, extend $f$ as in Theorem \ref{thm_main_param}.

Assume from now on that $S$ and $S_f$ are 1-cells over $Y$, with center $c$ and $c'$, and coset $\xi Q_{m,n}$ and $\xi' Q_{m',n'}$, respectively. Let $g$ be as in Theorem \ref{thm_preparation}, in particular $f$ and $g$ are equicompatible with $S$, and $(g(x,y)-c'(y))^b = e(y)(x-c(y))^a$ for every $(x,y)\in S$. 

Fix $y\in Y$ and let $B_{l,c(y),m,\xi}$ be a ball of $S$ above $y$. By Theorem \ref{thm_preparation} we can write $f_y(B_{l,c(y),m,\xi})=B_{l',c'(y),m',\xi'}=g_y(B_{l,c(y),m,\xi})$, where $B_{l',c',m',\xi'}$ is a ball of $S_f$ above $y$. Also, we have that $\ord(\partial f/\partial x) = \ord(\partial g/\partial x)$. Let $q=a/b$, then there are three different cases to consider, depending on whether $q=1$, $q<1$ or $q>1$.
\begin{description}
	\item[Case 1: $q=1$.] From equation \eqref{eq_formula_order} we have $\ord(\partial f(x,y)/\partial x) = \ord(e(y))$ for all $(x,y)\in S$. So for $x\in B_{l,c(y),m,\xi}$ we have
		\begin{align*}
			l' &= \ord(e(y)(x-c(y))+c'(y)-c'(y))\\
			&=\ord(e(y))+\ord(x-c(y))\\
			&=\ord(\partial f(x,y)/\partial x) + l,
		\end{align*}
			which implies $l'\geq l$, since $f$ is 1-Lipschitz in the first variable. In particular note that in this case $m=m'$, by Lemma \ref{lemma_l'm'ordf'lm}. This allows us to use the same extension as described in Remark \ref{rem_phi}, namely $\tilde{f}:K\times Y\to K: (x,y)\mapsto \hat{f}(\varphi_y(x),y)$, where $\hat{f}$ is as in the proof of Theorem \ref{thm_main_param} in the case that $S$ and $S_f$ are $1$-cells over $Y$, and $\varphi_y$ is as in Remark \ref{rem_phi} (again, remark that in Theorem \ref{thm_main_param} this extension is denoted with $\tilde{f}$). We prove that $\tilde{f}_y$ is 1-Lipschitz. Let $t_1\in\cup_l D_{l,c(y)}\cap X$ and $t_2\not\in\cup_l D_{l,c(y)}\cap X$, where $D_{l,c(y)} = \{x\in K\mid \ord(x-c(y))=l\}$. Then
		\begin{align*} 
			\abs{\tilde{f}_y(t_1)-\tilde{f}_y(t_2)} &= \abs{f(t_1,y)-c'(y)}\\
			&\leq \abs{t_1-c(y)}\\
			&\leq \abs{t_1-t_2},
		\end{align*}
		where the first inequality follows from $l'\geq l$ and the second from the non-Archimedean property.
	\item[Case 2: $q>1$.] Because $f$ is 1-Lipschitz in the first variable, we have $ \ord(\partial f/\partial x) \geq 0$, and together with \eqref{eq_formula_order} this gives the following lower bound: 
	\[l\geq -\ord(e(y)^{1/b}q)/(q-1).\]
	 Recall that $\ord(e(y)^{1/b}q)$ is short for $\ord(e(y))/b+\ord(q)$. On the other hand, as soon as $l\geq (m'-m-\ord(e(y)^{1/b}q))/(q-1)$, we have $l'\geq l$. Indeed, this follows immediately from Lemma \ref{lemma_l'm'ordf'lm} and from \eqref{eq_formula_order}. So up to partitioning $S$ into two cells over $Y$, we may assume that either $l'\geq l$ for all balls of $S$ above $y$, for every $y\in Y$, or that $S$ has at most $N$ balls above $y$, for every $y\in Y$, where $N$ does not depend on $y$. In the former case we can extend $f$ as we did in \textbf{Case 1}. In the latter case we can, after partitioning $Y$ in a finite number of definable sets, assume that there are \emph{exactly} $N$ balls of $S$ above $y$, for every $y\in Y$. Using (the remark after) Lemma \ref{lemma_glue} we may assume that there is exactly one ball of $S$ above $y$, for every $y\in Y$. By definable selection (see \cite{denef-vdd} and \cite{Dries1984-DRIATW}) there is a definable function $h:Y\to K$ such that for each $(x,y)\in S$ with $x\in K$ and $y\in Y$, $(h(y),y)\in S$. We then extend $f$ as follows:
		\[\tilde{f}:K\times Y\to K: (x,y)\mapsto\begin{cases}f(x,y) &\text{if }(x,y)\in S,\\ f(h(y),y)&\text{if }(x,y)\not\in S.\end{cases}\]
		Fix $y\in Y$, we show that $\tilde{f}_y$ is 1-Lipschitz. Recall that by the argument given above, $S_y$ is a ball in $K$. The only nontrivial case to consider is the following. Let $t_1\in S_y$ and $t_2\not\in S_y$, then
		\begin{align*}
			\abs{\tilde{f}_y(t_1)-\tilde{f}_y(t_2)} &= \abs{f(t_1,y)-f(h(y),y)}\\
			&\leq \abs{t_1-h(y)}\\
			&<\abs{t_1-t_2},
		\end{align*}
		where the last inequality holds because of the non-Archimedean property and the fact that $t_1$ and $h(y)$ are both contained in the ball $S_y$, and $t_2$ is not.
	\item[Case 3: $q<1$.] This case is similar to \textbf{Case 2}, where now one finds an upper bound for $l$ instead of a lower bound. The proof is omitted. \qedhere
\end{description}
	\end{proof}
	\begin{remark}
	Note that we proved the main theorem for semi-algebraic and subanalytic structures on $K$. It is, for now, unclear whether the main theorem could also hold in other structures on $K$, such as, for example, $P$-minimal structures, as defined by Haskell and Macpherson in \cite{has-mac-97}. Also, it is unclear whether the extension that we constructed could be used to extend a definable function $f:X\subset K^r\to K^s$ that is $\lambda$-Lipschitz in \emph{all} variables to a definable function $\tilde{f}:K^r\to K^s$ that is $\lambda$-Lipschitz in all variables. For now, there is no evidence towards either a positive or a negative answer to this question.
	\end{remark}

\bibliographystyle{plain}
\bibliography{references}
\vspace*{10pt}
\textsc{Tristan Kuijpers\\ KU Leuven, Department of Mathematics, Celestijnenlaan 200B, 3001 Leuven, Belgium}\\
\textit{E-mail:} tristan.kuijpers@wis.kuleuven.be
\end{document}